\documentclass[12pt]{amsart}
\newif\ifpdf
    \ifx\pdfoutput\undefined
    \pdffalse 
    \else
    \pdfoutput=1 
    \pdftrue
    \fi

    \ifpdf
    \usepackage[pdftex]{graphicx}
    \else
    \usepackage{graphicx}
    \fi
\usepackage{amsmath,amssymb,amsthm,epsfig}

\newcommand\g{\gamma}

\newcommand\F{\mathcal F}

\newcommand\C{\mathcal C}

\newcommand\x{x_0}

\newcommand\inv{x_0^{-1}}

\newcommand\Ll{L_L}

\newcommand\rni{R_{NI}}

\newtheorem{theorem}{Theorem}[section]

\newtheorem{lemma}[theorem]{Lemma}

\newtheorem{corollary}[theorem]{Corollary}

\input{epsf.tex}

\begin{document}
\ifpdf
    \DeclareGraphicsExtensions{.pdf, .jpg, .tif}
    \else
    \DeclareGraphicsExtensions{.eps, .jpg}
    \fi

\title{Thompson's group $F$ is not almost convex}
\author{Sean Cleary and Jennifer Taback}
\thanks{The first author acknowledges support from PSC-CUNY grant \#63438-0032}
\thanks{The second author would like to thank the University of
Utah for their hospitality during the writing of this paper.}

\begin{abstract}
We show that Thompson's group $F$ does not satisfy Cannon's almost
convexity condition $AC(n)$ for any integer $n$ in the standard
finite two generator presentation.  To accomplish this, we construct
 a family of pairs of elements at distance $n$ from the identity
and distance 2 from each other, which are not connected by a path lying
inside the $n$-ball of length less than $k$ for increasingly large $k$.
Our techniques rely upon Fordham's method for calculating the length
of a word in $F$ and upon an analysis of the generators' geometric
actions on the tree pair diagrams representing elements of $F$.

\end{abstract}

\maketitle

\section{Introduction}

\label{sec:intro}

Cannon \cite{cannon} introduced the notion of almost convexity for
a group $G$ with respect to a finite generating set $X$. This
finite generating set $X$ determines a word metric $d_X$ for $G$
and its Cayley graph. $G$ is {\em almost convex $(k)$} or  {\em
$AC(k)$} with respect to $X$ if there is a number $N(k)$ so that
for all positive integers $n$, given two elements $y$ and $z$ in
the ball $B(n)$ of radius $n$ with $d_X(y,z) \leq k$,
there is a path $\gamma$ from $y$ to $z$ of length at most
$N(k)$ which lies entirely in $B(n)$.   Cannon showed  that if a
group $G$ is $AC(2)$ with respect to a finite generating set then
$G$ is $AC(k)$ for $k \geq 2$ and thus a group satisfying $AC(2)$
is called {\em almost convex}.  Almost convexity allows
algorithmic construction of $B(n+1)$ from $B(n)$ by making it
sufficient to consider only a finite set of possible ways that an
element  in $B(n+1)$ can be obtained from different elements of
$B(n)$.

A number of families of groups have been shown to be
almost convex.  Cannon \cite{cannon} showed that 
 hyperbolic groups are almost convex and
 that amalgamated products of almost convex groups are almost convex.
Stein and Shapiro \cite{steinshapiro} showed that 
 fundamental groups of closed three manifolds whose
geometry is not modelled on $Sol$  are almost convex.
Other families of groups have been shown not to be almost convex.
Cannon, Floyd, Grayson and Thurston \cite{solcannon}
showed that  fundamental groups of manifolds
with $Sol$ geometry \cite{solcannon} are not almost convex,
and Miller and Shapiro showed that the 
 solvable Baumslag-Solitar groups
$BS(1,n)$ \cite{cfmshapiro} are not almost convex.
Unfortunately, the property of almost convexity can depend
upon presentation.   Thiel \cite{thiel} showed that generalized
Heisenberg groups are not almost convex with respect to
their standard presentations, but are almost convex with
respect to some alternate presentations.

Although Thompson's group $F$ has been studied extensively in
many branches of mathematics, the metric properties of $F$ were
poorly-understood until recently.  Burillo \cite{burillo} and
Burrilo, Cleary and Stein \cite{bcs} developed estimates for
measuring distance in $F$, and Fordham \cite{blake} developed a
remarkable method for computing  distance in $F$. 

We prove below that $F$ does not satisfy
Cannon's $AC(2)$ property in its standard finite generating set,
and thus is not almost convex with respect to that generating set.

Thompson's group $F$ has a number of different manifestations.
Originally discovered in logic as the group of automorphisms of a
free algebra by Thompson \cite{thomp}, $F$ also has connections with homotopy
theory developed by Freyd and Heller \cite{fh1,fh2},  groups of homeomorphisms of the
interval studied by Brin and Squier \cite{bs:pl}  and Brown and Geoghegan
\cite{bg:thomp} and diagram groups defined by Guba and Sapir \cite {diag}. 
Cannon, Floyd and Parry \cite{cfp} give an introduction and summarize many of the
remarkable properties of $F$.

Thompson's group $F$ has the infinite presentation  ${\mathcal P}$
given by
$$
{\mathcal P} = \langle x_k, \ k \geq 0 | x_i^{-1}x_jx_i = x_{j+1}
\ \text{ if }i<j \rangle.$$  We can see that the lower index
generators conjugate the higher-index generators by incrementing
their indices.  Since $x_0$ conjugates $x_1$ to $x_2$ and successively
to all higher index generators, it is clear
that $F$ is finitely generated.  In fact, all of the infinitely
many relators in ${\mathcal P}$ are consequences of a basic set of two relators.
Thus, there is the following standard finite presentation ${\mathcal F}$ for
$F$:
$$
{\mathcal F} = \langle x_0,x_1 |
[x_0x_1^{-1},x_0^{-1}x_1x_0],[x_0x_1^{-1},x_0^{-2}x_1x_0^2]
\rangle.$$

We  prove the following theorem:
\begin{theorem}
\label{thm:notAC} Thompson's group $F$ does not satisfy Cannon's
almost convexity condition $AC(2)$ in the finite presentation
${\mathcal F}$.
\end{theorem}

We immediately obtain the corollary:
\begin{corollary}
\label{cor:notACn} Thompson's group $F$ does not satisfy Cannon's
almost convexity condition $AC(n)$ for any positive integer $n >
2$ in the finite presentation $\F$.

\end{corollary}

\section{Background on $F$}

\label{sec:background}

Analytically, we define $F$ as the group of orientation-preserving
piecewise-linear homeomorphisms from $[0,1]$ to itself where each
homeomorphism has only finitely many singularities of slope, all
such singularities lie in the dyadic rationals ${\bf Z}[\frac12]$,
and, away from the singularities, the slopes are powers of $2$.

Combinatorially, $F$ has the infinite and finite presentations
given above. There is a convenient set of normal forms for
elements of $F$ in the infinite presentation ${\mathcal P}$ given
by $x_{i_1}^{r_1} x_{i_2}^{r_2}\ldots x_{i_k}^{r_k} x_{j_l}^{-s_l}
\ldots x_{j_2}^{-s_2} x_{j_1}^{-s_1} $ with $r_i, s_i >0$,
$i_1<i_2 \ldots < i_k$ and $j_1<j_2 \ldots < j_l$. This normal
form is unique if we further require that when both $x_i$ and
$x_i^{-1}$ occur, so does $x_{i+1}$ or  $x_{i+1}^{-1}$, as
discussed by Brown and Geoghegan \cite{bg:thomp}.  In what follows, when we refer to a
word in normal form, we always mean the unique normal form.

The geometric description of $F$ is in terms of tree pair
diagrams. A tree pair diagram is a pair of rooted binary trees
with the same number of leaves, as described in  \cite{cfp}.   We
number the leaves of each tree from left to right, beginning with
0. We refer to an interior node together with the two downward-directed
edges from the node as a {\em caret}.  We define the right (respectively left) child of a caret
$C$ to be the caret  $C_R$ (respectively $C_L$) which is attached to the right (left)
downward edge of caret $C$.

Each tree in a tree pair can
be regarded as a set of instructions for successive subdivision of
the unit interval: the root caret subdivides the interval in half,
a right child of the root subdivides $[ \frac{1}{2},1]$ in half,
and so on.  This gives a correspondence between elements of $F$ in
the geometric description and the analytic description as follows.
Let $(T_-,T_+)$ be a pair of trees each with $n$ leaves.
Each tree determines a subdivision of $[0,1]$
into $n$ subintervals. The
tree pair  $(T_-,T_+)$ corresponds to the piecewise linear homomorphism which
maps the  subintervals of the $T_-$ subdivision to the subintervals of
the $T_+$ subdivision, in order.  This equivalence and the group
operation are described in \cite{cfp}.  We refer to $T_-$ as the
{\em negative tree} and $T_+$ as the {\em positive tree} of the
pair $(T_-,T_+)$.

A tree pair diagram is {\em unreduced} if each of $T_-$ and $T_+$
contain a caret with leaves numbered $m$ and $m+1$, and it
is {\em reduced} otherwise.
Note that there are many tree pair diagrams representing the same
element of $F$ but there is a unique reduced tree pair diagram for
each element of $F$.  When we write $(T_-,T_+)$ to represent an element of $F$,
we are assuming that the tree pair is reduced.

If $x = (T_-,T_+)$ is a reduced pair of trees representing $x$,
the normal form for $x$ can be constructed by the following
process, described in \cite{cfp}. Beginning with the tree pair $(T_-,T_+)$, we number the
leaves of $T_-$ and $T_+$ from left to right, beginning with $0$.
The {\em exponent} of the leaf labelled $n$, written $E(n)$, is
defined as the length of the maximal path consisting entirely of
left edges from $n$ which does not reach the right side of the
tree. Note that $E(n)=0$ for a leaf labelled $n$ which is a right child of a
caret, as there is no path consisting entirely of left edges
originating from $n$.

We compute $E(n)$ for all  leaves in $T_-$, numbered $0$ through
$m$. The negative part of the normal form for $x$ is then
$x_m^{-E(m)} x_{m-1}^{-E(m-1)} \cdots x_1^{-E(1)} x_0^{-E(0)}$. We
compute the exponents for the leaves of the positive tree and thus
obtain the positive part of the normal form as
$x_0^{E(0)} x_{1}^{E(1)} \cdots x_m^{E(m)}$.
 Many of the exponents
may be $0$, and after deleting these, we can index the remaining
terms to correspond to the normal form given above, as detailed in
\cite{cfp}.

\begin{figure}

\includegraphics[width=5in]{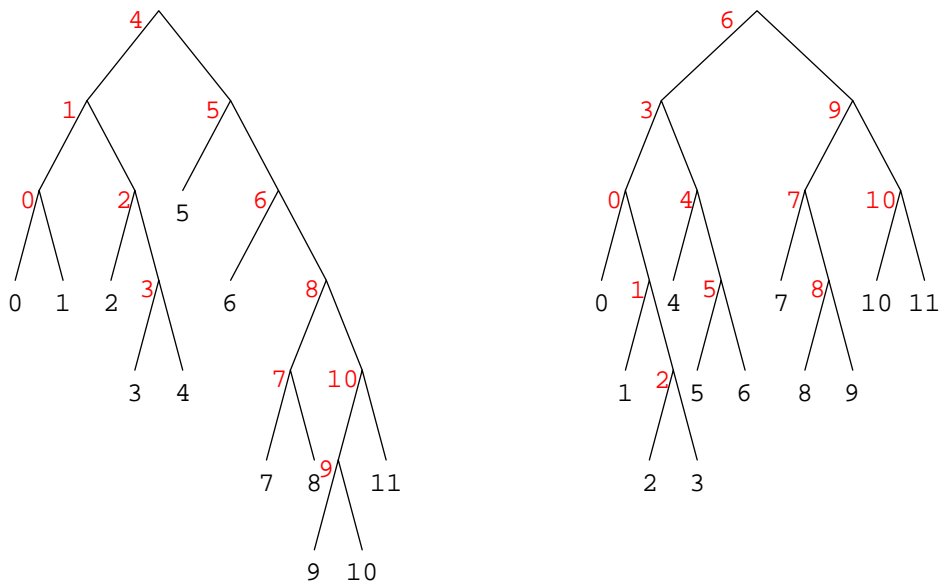}\\
\caption{Tree pair diagram for
 $x_0^2 x_1 x_2 x_4 x_5 x_7 x_8 x_9^{-1} x_7^{-1} x_3^{-1} x_2^{-1} x_0^{-2}$
with carets and leaves numbered. \label{treepairexample}}
\end{figure}

In the tree pair diagram in Figure \ref{treepairexample}, the exponent $E(0)$ of the
leaf labelled $0$ of $T_-$ is 2 since there is a path of two left
edges from leaf 0 which does not reach the right hand side of the
tree. The third left edge emanating from leaf 0 touches the
right-hand side of the tree and thus does not contribute to the
exponent. The exponents of all the leaves of $T_-$ are, in order,
$2,0,1,1,0,0,0,1,0,1,0,0$, and the exponents of the leaves of
$T_+$ are, in order, $2,1,1,0,1,1,0,1,1,0,0,0$.  Using these
exponents, and omitting any which are $0$, we see that the tree pair
diagram of Figure \ref{treepairexample} represents the word $x_0^2 x_1 x_2 x_4 x_5 x_7
x_8 x_9^{-1} x_7^{-1} x_3^{-1} x_2^{-1} x_0^{-2}$, in normal form.

If $R$ is a caret on the right side side of the tree  with a single left leaf labelled $k$, then
$E(k) = 0$ by definition.  Thus, right carets with no left subtrees
can be added to the rightmost carets of  either $T_-$ or $T_+$
without affecting the normal
form to ensure that both trees have the same number of carets.

Similarly, given an element $x$ in normal form with respect to the
infinite generating set, it is possible to construct a tree pair
diagram  ($T_-, T_+$) so that each leaf has the correct exponent.
In particular, the number of left edges of $T_-$ emanating from the root caret
is one more than the exponent of $x_0^{-1}$ in the normal form and
the number of left edges of $T_+$ emanating from the root caret is one more than
the exponent of $x_0$ in the normal form for $x$.

The processes described above relate the normal form of words in
$F$ in the infinite presentation ${\mathcal P}$ to the tree pair
representation. For many questions involving the geometry of $F$,
we must consider the length of words in $F$ with respect to a
metric arising from a finite generating set.  Burillo
\cite{burillo} presented a way of estimating the word length
$| x |_{\F}$ in the finite generating set $\F$ from the normal
form, which was refined by Burillo, Cleary and Stein in \cite{bcs}.

\begin{theorem}[Burillo \cite{burillo}, Prop.\ 2 ; Burillo, Cleary and Stein \cite{bcs}, Theorem 1]
\label{thm:burillo} Let $w \in F$ have normal form $w =
x_{i_1}^{r_1}\cdots x_{i_n}^{r_n}x_{j_m}^{-s_m} \cdots
  x_{j_1}^{-s_1}$, and let $D(w) = r_1 + r_2 + \cdots + r_n + s_1
  + s_2 + \cdots + s_m + i_n + j_m$.  Then $$ \frac{D(w)}{3} \leq
  |x|_{\F} \leq 3D(w).$$
\end{theorem}
 Burillo, Cleary and Stein \cite{bcs} also estimated of
the length $|x|_{\F}$ of a word $x$ given by a tree pair diagram
in terms of the number of carets $N(x)$ in either tree.

\subsection{Fordham's method of calculating word length}
 Fordham \cite{blake} presents a method of calculating the
exact word length in $F$ given a reduced pair of trees
representing an element $x \in F$.  We make some preliminary
definitions before explaining Fordham's technique.

Let $T$ be a finite rooted binary tree.  The {\em left side} of
$T$ is the maximal path of left edges beginning at the root of
$T$. Similarly, we have the {\em right side} of $T$.  A caret in $T$
is a {\em left caret} if its left edge is on the left side of the
tree, a {\em right caret} if it is not the root and its right edge
is on the right side of the tree, and an {\em interior caret}
otherwise. The carets in $T$ are numbered according to the
 infix ordering of nodes. We begin numbering with leaf $0$ as the
leftmost leaf and caret $0$ the left caret whose left child is
leaf $0$. We number the left children of a caret before the caret
itself, and number the right children after numbering the caret.  The
trees in Figure \ref{treepairexample} have their carets numbered according to this
method.

Fordham classifies carets into seven disjoint types:
\begin{enumerate}

\item
$L_0$.  The first caret on the left side of the tree, with
caret number $0$. Every tree has exactly one caret of type $L_0$.

\item
$\Ll$.  Any left caret other than the one numbered $0$.

\item
$I_0$. An interior caret which has no right child.

\item
$I_R$. An interior caret which has a right child.

\item
$R_I$. Any right caret numbered $k$ with the property that caret
$k+1$ is an interior caret.

\item
$\rni$. A right caret which is not an $R_I$ but for which there is a
higher numbered interior caret.

\item
$R_0$. A right caret with no higher-numbered interior carets.

\end{enumerate}

 The root caret is always considered to be a left caret of type $\Ll$ unless
it has no left children, in which case it is the $L_0$ caret.

Working from caret 0 to caret 10, in infix order, in the tree $T_-$ from Figure \ref{treepairexample}, we
see that the carets are of types $$L_0, \Ll, I_R, I_0, \Ll, \rni, R_I,
I_0, R_I, I_0 \text{ and } R_0.$$  The carets in the tree $T_+$ of
Figure \ref{treepairexample}, in infix order, are of types $$L_0,
I_R, I_0, L_L,I_R,I_0,L_L,I_R,I_0,R_0 \text{ and } R_0.$$

The main result of  Fordham \cite{blake} is that the word length $|x|_{F}$
of $x = (T_-,T_+)$ can be computed from knowing the caret types of
the carets in the two trees, as long as they form a reduced pair,
 via the following process.  We number
the $k+1$ carets according to the infix method described above,
and for each $i$ with $0 \leq i \leq k$ we form the pair of caret
types consisting of the type of caret number $i$ in $T_-$ and the
type of caret number $i$ in $T_+$. The single caret of type $L_0$
in  $T_-$ will be paired with the single caret of type $L_0$ in
$T_+$, and for that pairing we assign a weight of 0. For all other
caret pairings, we assign weights according to the following symmetric
table:

\begin{center}

\begin{tabular}{|c|c|c|c|c|c|c|}

\hline
 & $R_0$ & $\rni$ & $R_I$ & $\Ll$ & $I_0$ & $I_R$ \\
 \hline

 $R_0$ & 0 & 2 & 2 & 1 & 1 & 3 \\ \hline
 $\rni$ & 2 & 2 & 2 & 1 & 1 & 3 \\ \hline
 $R_I$ & 2 & 2 & 2 & 1 & 3 & 3 \\ \hline
 $\Ll$ & 1 & 1 & 1 & 2 & 2 & 2 \\ \hline
 $I_0$ & 1 & 1 & 3 & 2 & 2 & 4 \\ \hline
$I_R$ & 3 & 3 & 3 & 2 & 4 & 4 \\ \hline
\end{tabular}

\end{center}
Fordham's remarkable result is that the sum of these weights is exactly the
 length of the word in the word metric arising from the finite generating set.

\begin{theorem}[Fordham \cite{blake}, Theorem 2.5.1]
\label{thm:blake} Given a word $w \in F$ described by the reduced tree
pair diagram $(T_-,T_+)$, the length $|w|_{\F}$ of the word with respect
to the generating set $\F$ is the sum of
the weights of the caret pairings in $(T_-,T_+)$.
\end{theorem}

Considering the word $w$ in Figure \ref{treepairexample}, we see that the carets numbered
zero have type pairing $(L_0,L_0)$, which has weight $0$.  The
carets numbered $1$ have types $(L_L,I_R)$ which contributes 2 to
the weight of the word.  The total weight of the word is easily
computed to be 0+2+4+2+2+1+1+4+3+1+0=20.  Thus, the length of $w$ in the
word metric $| w |_{\F}$ is $20$.

The proofs in \S \ref{sec:proofs} rely heavily on this technique of
Fordham.  Namely, we use the fact that we can apply a generator to
a given word, whose length we know, and the change in caret types,
which is easily seen, exactly determines the change in word
length.

\subsection{Action of the generators on an element of $F$}

\label{sec:action}

We begin with a lemma from Fordham \cite{blake} which states under fairly broad conditions,
that when applying a generator to a tree pair $(T_-,T_+)$
exactly one pair of caret types will change.  In \S
\ref{sec:family}, we construct a special family of elements which
will provide the counterexamples to almost convexity for $F$.
These elements are constructed to satisfy the conditions of the lemma below.

\begin{lemma}[Fordham \cite{blake}, Lemma 2.3.1]
\label{lemma:conditions}
Let $(T_-,T_+)$ be a reduced pair of trees, each having $m+1$
carets, representing an element $x \in \F$, and $\alpha$ any
generator of $\F$.
\begin{enumerate}
\item
If $\alpha = \x$, we require that the left subtree of the root of $T_-$
is nonempty.

\item
If $\alpha = \inv$, we require that the right subtree of the root of $T_-$
is nonempty.

\item
If $\alpha = x_1$, we require that the left subtree of the right child
of the root of $T_-$ is nonempty.

\item
If $\alpha = x_1^{-1}$, we require that the right subtree of the right child
of the root of $T_-$ is nonempty.

\end{enumerate}
If the reduced tree pair diagram for $x \alpha$ also
has $m+1$ carets,
then there is exactly one $i$ with $0 \leq i \leq
m$ so that the pair of caret types of caret $i$ changes when
$\alpha$ is applied to $x$.
\end{lemma}

We now begin to understand geometrically the action of a generator
of $\F$ on a reduced tree pair $(T_-,T_+)$, and the corresponding change in
normal form. We will generally  assume that the conditions of lemma
\ref{lemma:conditions} are met by the generic elements with which we begin.

Let $C_R$ denote the caret which is the right child of the root
caret $R$ of $T_-$, and $C_{RR}$ and $C_{RL}$ the right and left
carets, respectively, of $C_R$. Similarly, let $C_L$ denote the
left child of the root caret of $T_-$, and $C_{LL}$ and $C_{LR}$
its left and right children.  Figures $2,3$ and $4$ will be useful
in understanding the geometric interpretation of the action of the
generators on an element of $F$.  In all of these figures, the
letters $a,b$ and $c$  represent (possibly empty) subtrees of the
given tree.

We first understand the action of the generator $\inv$ on a tree
pair $(T_-,T_+)$ representing an element $w \in F$.  Consider $w$
written in normal form as $w = x_{i_1}^{r_1}\cdots
x_{i_n}^{r_n}x_{j_m}^{-s_m} \cdots x_{j_1}^{-s_1}$.  Then the
element $w \inv$ is still in normal form (unless we are in the
degenerate case where $x=x_0^m$.)  Recall from \S
\ref{sec:background} that the exponent of $\inv$ in the normal
form is one less than the number of left edges of the tree $T_-$.
Thus, increasing the exponent of $\inv$ by $1$ adds a left edge to
$T_-$.

The numbering of the leaves and carets after this new edge is added must
remain the same, since the normal form (and hence the exponents of
the leaves) changes in a single place.  Thus, with the extra edge
in $T_-$, $C_R$ becomes the new root caret.  The  left
subtree $C_{RL}$ of  $C_R$, which contains carets with smaller numbers than
$C_R$, must become the right subtree of the old root caret, which is
now at position  formerly occupied by $C_L$. The left caret $C_L$ is moved
down and to the left and remains a left caret, now in the position formerly
occupied by $C_{LL}$ and so on.
This tree transformation is also called a {\em counterclockwise rotation} 
or {\em left rotation} based at the
root.  Figure 2 shows the negative trees $T_-$ for the elements
$w$ and $w \inv$ and illustrates a counterclockwise rotation based
at the root.

\begin{figure}\includegraphics[width=3in]{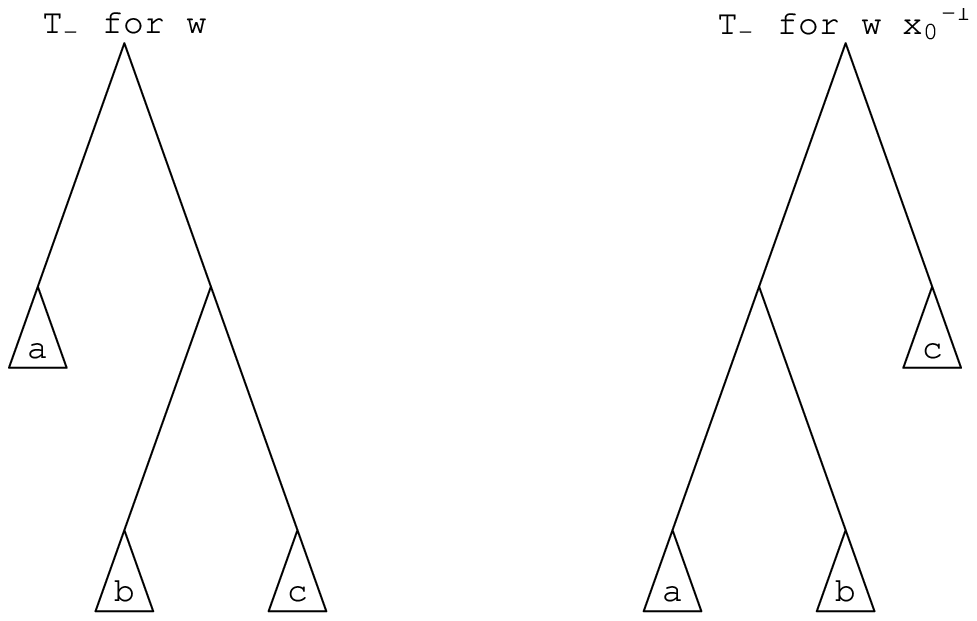}\\
\caption{Rotation at the root induced by applying $\inv$ to $T_-$.
\label{x0invaction}}
\end{figure}

\begin{figure}\includegraphics[width=3in]{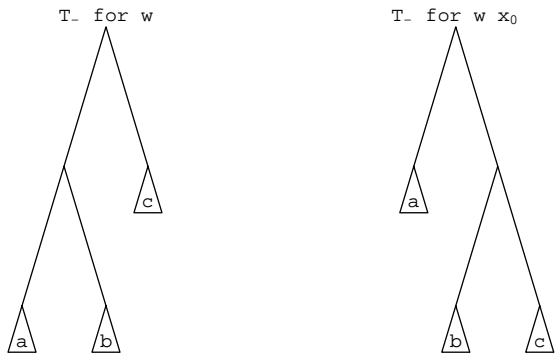}\\
\caption{Rotation at the root induced by applying $\x$ to $T_-$. \label{x0action}}
\end{figure}

When we consider the action of $\x$ on $w=(T_-,T_+)$, we can
assume, according to lemma \ref{lemma:conditions}, that $T_-$
has at least two left edges, equivalently, that the exponent of
$\inv$ in the normal form of $w$ is at least $1$. Applying the
generator $\x$ cancels one $\inv$ in the normal form. This
corresponds to the tree $T_-$ losing a left edge, and thus the
caret  $C_L$ becomes the root caret and the former root
caret $R$ moves to the position of $C_R$.  The initial right
subtree $C_{LR}$ of $C_L$ becomes the left subtree of $R$ in order to
preserve the numbering of the carets.  This is a clockwise (or right)
rotation based at the root of $T_-$ and is illustrated in Figure
\ref{x0action} .

It is more difficult to visually understand the action of $x_1$
and $x_1^{-1}$ on the pair $(T_-,T_+)$ corresponding to $w$, as it
is more difficult to see how these generators change the normal
form. Using the terminology given above, the following lemmas show
that the generators $x_1$ and $x_1^{-1}$ perform counterclockwise and
clockwise rotations around the node $C_R$.

We begin with a lemma relating the action of $x_1^{-1}$ on
$(T_-,T_+)$ to the normal form of the corresponding element $w \in
\F$.

\begin{lemma}[The normal form of $wx_1^{-1}$]
\label{lemma:leafnumber1} Let $w  \in F$ be represented by the
tree pair $(T_-,T_+)$, and have normal form $x_1^{r_1} \cdots
x_{i_n}^{r_n}x_{j_m}^{-s_m} \cdots x_{j_1}^{-s_1}$.  Then $w
x_1^{-1}$ has normal form
\begin{equation}
\label{eqn:nf1}
 x_{i_1}^{r_1} \cdots x_{i_n}^{r_n}x_{j_m}^{-s_m}
\cdots x_{j_{q+1}}^{-s_{q+1}} x_{\alpha}^{-1} x_{j_q}^{-s_q}
\cdots x_{j_1}^{-s_1},
\end{equation}
where we might have $\alpha  =  j_{q+1}$.  If the root of $T_-$
has right and left subtrees $T_R$ and $T_L$ respectively, then
$\alpha$ is smallest leaf number in $T_R$.
\end{lemma}

\begin{proof}
We consider the proof in two cases.  In the first case, if $j_1
\neq 0$ then $\alpha = 1$ and the expression $x_1^{r_1} \cdots
x_{i_n}^{r_n}x_{j_m}^{-s_m} \cdots x_{j_1}^{-s_1} x_1^{-1}$ is in
normal form.  In this case, $T_-$ has a single left edge, with
leaf labelled $0$, and the first left leaf of the first right
subtree will be labelled $1$.

In the second case we assume that $j_1 = 0$.  Then the relators in
${\mathcal P}$ imply that $\alpha = 1+ s_1 + s_2 + \cdots + s_l$,
where $l$ is the first index satisfying $j_{l+1} \geq 1+ s_1 + s_2 +
\cdots + s_l$.  It remains to show that this is the label of the
leftmost leaf of the first right subtree of $T_-$.

Let $T_L$ and $T_R$ be the left and right subtrees of the root
caret of $T_-$.  We consider the number of interior carets in
$T_L$. If $T_L$ is empty, then we are in the first case discussed
above.

If $T_L$ has no interior carets, but is not empty, then the number
of left edges in $T_L$ is $n$, for some $n$, and thus the last
leaf number in $T_L$ is $n$ as well.  So the first leaf number in
$T_R$ is $n+1$. Given this form of $T_-$, we see that the normal
form of $x$ must end with $x_{j_2}^{-s_2} \x^{-n}$ where $j_2 \geq
n+1$. Thus, using the relators to put $x_1^{-1}$ into its proper
position in the normal form, we see that it becomes
$x_{1+n}^{-1}$, agreeing with the statement of the lemma.

If $T_L$ has a single interior caret, then the total number of
left edges of $T_L$ is $n+1$, where $n$ again represents the
length of the left side of $T_L$.  The interior caret also adds an
additional leaf, and thus the highest numbered leaf of $T_L$ is
$n+1$. We know that $x_1^{-1}$ becomes $x_{n+1}^{-1}$ when it is
moved left past the $x_0^{-n}$.  However, $n+1$ is now the highest
numbered leaf in $T_L$.  The extra left leaf added by the single
interior caret corresponds to a letter in the normal form of $x$
whose index is smaller than $n+1$, thus when the $x_{n+1}^{-1}$ is
moved left past this letter, it becomes an $x_{n+2}^{-1}$.  Since
there are no other interior carets in $T_L$, the next possible
index of a letter in the normal form of $x$ is $n+2$.  Thus
$x_{n+2}^{-1}$ is now in place in the normal form, so $\alpha =
n+2$, and $n+2$ is the first leaf number of $T_R$, as required.

If $T_L$ has two interior carets, then there are $n+2$ left edges
in $T_L$ and the highest leaf number in $T_L$ is $n+2$.  Moving left
past $\x^{-n}$, the $x_1^{-1}$ first becomes $x_{n+1}^{-1}$, as in
the previous case.  Again, we see that $n+1$ is a leaf number in
$T_L$.  Then, since there is a single leaf numbered higher than
$n+1$, the are not enough leaves to have the remaining two carets
have leaves numbered higher than $n+1$.  So the first interior
caret must have a leaf with a lower number than $n+1$,
corresponding to a letter in the normal form of $x$ with index
smaller than $n+1$. Thus $x_{n+1}^{-1}$ must be moved left past this
element as well, making it $x_{n+2}^{-1}$.  Now, $n+2$ is the
highest leaf number in $T_L$, so the second interior caret must
again appear before leaf number $n+2$; that is, it corresponds to a
letter in the normal form of $x$ with index smaller than $n+2$.
Moving the $x_{n+2}^{-1}$ past left this letter, we get
$x_{n+3}^{-1}$.  Since there are no more interior carets in $T_L$,
there are no other letters in the normal form with index less than
$n+3$, so we must have $x_{n+3}^{-1}$ in its place in the normal
form.  Again, we see that $n+3$ is the first leaf number in $T_R$.

In summary, each additional interior caret adds a letter to the
normal form with smaller index than $n+1$; thus the $x_1^{-1}$
must be moved left past these letters to obtain the normal form.  We
can continue this method to apply to an arbitrary number of
interior carets in $T_L$, proving the lemma.
\end{proof}

\begin{lemma}[The normal form of $wx_1$]
\label{lemma:leafnumber2} Let $w$ satisfy the conditions of lemma
\ref{lemma:conditions} and have normal form $x_1^{r_1} \cdots
x_{i_n}^{r_n}x_{j_m}^{-s_m} \cdots x_{j_1}^{-s_1}$.  Then $w x_1$
has normal form:
\begin{equation}
\label{eqn:nf2} x_1^{r_1} \cdots x_{i_n}^{r_n}x_{j_m}^{-s_m}
\cdots x_{j_l}^{-(s_l-1)} \cdots x_{j_1}^{-s_1},
\end{equation}
for some index $j_l$, in which case $j_l$ is the smallest
leaf number in the right subtree of $T_-$.

\end{lemma}

\begin{proof}
As in the proof of lemma \ref{lemma:leafnumber1}, we use the
relators of ${\mathcal P}$ to move $x_1$ to the generator
$x_{\alpha}$ where $\alpha = 1+ s_1 + s_2 + \cdots + s_l$, and $l$
is the first index satisfying the inequality $j_{l+1} \geq 1+ s_1 +
s_2 + \cdots + s_l$, or $\alpha = 1+s_1 +s_2 + \cdots + s_m$. From
the proof of lemma \ref{lemma:leafnumber1} we again know that
$\alpha$ is the number of the leftmost leaf of the first right
subtree of $T_-$.

According to lemma \ref{lemma:conditions}, the left subtree of
$C_L$ is nonempty, so there is a leaf labelled $\alpha$ with
exponent at least $2$, i.e. there is an index $j_k = \alpha$ in
the normal form of $w$.  Thus the exponent of $x_{j_k}$ decreases
by $1$ because the $x_{\alpha}$ cancels one $x_{j_k}^{-1}$ letter
giving the normal form (\ref{eqn:nf2}).

\end{proof}

\begin{lemma}[The action of $x_1^{-1}$ on $T_-$]
\label{lemma:x1inverse} The generator $x_1^{-1}$ when applied to
an element $w$ of $F$ represented by a tree pair $(T_-,T_+)$
which satisfies the conditions of lemma \ref{lemma:conditions} leaves $T_+$ unchanged,
and affects $T_-$ as follows: $C_{RR}$ becomes the right child of the root caret,
and $C_R$ becomes the left child of $C_{RR}$.  All other carets remain unchanged.
\end{lemma}

\begin{figure}\includegraphics[width=3in]{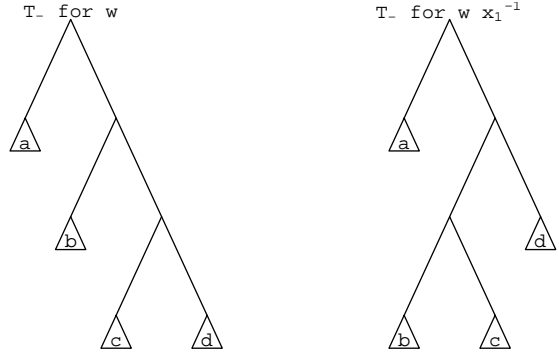}\\
\caption{Left rotation around $C_R$ induced by applying $x_1^{-1}$
\label{x1invaction}}
\end{figure}

\begin{proof}Let $\alpha$ be the number of the leftmost leaf in the right
subtree of the root of $T_-$.  It follows from lemma
\ref{lemma:leafnumber2} that the exponent of $x_{\alpha}$ in the
normal form of $x$ is increased by $1$; that is, the exponent
$E(\alpha)$ of the leaf $\alpha$ is increased by $1$, which means
there is one more left edge emanating from $C_R$ in $T_-$ and
terminating at $\alpha$.  Since the numbering of the carets is
preserved, because the normal form changes in a single letter, and
begins at the far left of the right subtree of the root caret, we
see that $C_R$ is now an interior caret.  To preserve the
numbering of the leaves and carets, the left subtree of $C_{RR}$
must become the right subtree of $C_R$, because these carets are
numbered higher than $C_R$ but lower than $C_{RR}$.  This leaves
$C_{RR}$ as the right child of the root caret. All remaining
subtrees are left unchanged.
\end{proof}

\begin{lemma}[The action of $x_1$ on $T_-$]
\label{lemma:x1} The generator $x_1$ when applied to an element $w
\in  F$ represented by a tree pair $(T_-,T_+)$ satisfying the
conditions of lemma \ref{lemma:conditions}  leaves $T_+$
unchanged, and in $T_-$, causes $C_{RL}$ to become the right child
of the root and $C_R$ to become the right child of $C_{RL}$. All
other carets remain unchanged.
\end{lemma}

\begin{figure}\includegraphics[width=3in]{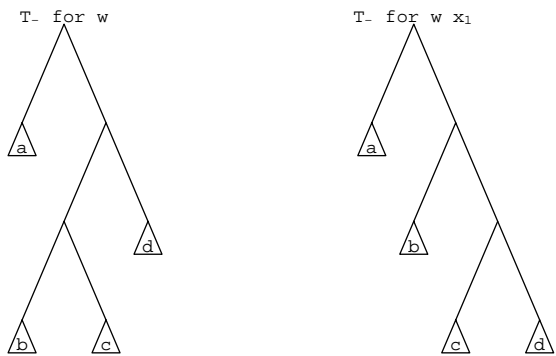}\\
\caption{Right rotation around $C_R$ induced by applying $x_1$
\label{x1action}}
\end{figure}

\begin{proof}
The normal form of $w x_1$ is of the form (\ref{eqn:nf2}) given in
lemma \ref{lemma:leafnumber2}.  From lemma \ref{lemma:leafnumber1}
we know that the index $j_l$ is the number of the leftmost leaf in
the left subtree of $C_R$ in $T_-$. From the change in normal form
we see that the exponent of $x_{j_l}$ decreases by $1$ and thus in
$T_-$ the exponent $E(j_l)$ decreases by $1$.  Thus, there is one
fewer left edge emanating from $C_R$ ending in the leaf numbered
$j_l$. Accordingly, the right subtree of $C_{RL}$ is moved to the
right side of $T_-$, without changing the numbering of the carets.
Thus $C_{RL}$ is now the right child of the root, and $C_R$ is the
left child of $C_{RL}$.

\end{proof}

Notice that in all of the descriptions above, the tree $T_+$ is
not affected by the action of a generator.
This is not true in general for reduced tree pair diagrams not satisfying the
conditions of lemma \ref{lemma:conditions}.
In general, $T_+$ can be
affected in exactly three ways:

\begin{enumerate}
\item when $T_-$ has a single left edge, and the generator is
$\x$,

\item
when the left  subtree of $C_{R}$ of $T_-$ is empty, and the generator
is $x_1$, or

\item
if the generator is $\alpha$ and the pair of trees corresponding
to $x \alpha$ is not reduced.
\end{enumerate}
We choose the family of words which will provide the
counterexamples to almost convexity so that the conditions of lemma
\ref{lemma:conditions} are always satisfied.

\section{A special family of elements}

\label{sec:family}

We define a family $\C(k)$, with integral $k \geq 2$, of elements of $F$
which we will use to prove that $F$ is not $AC(2)$, and thus not
$AC(n)$.  We first define what the negative tree $T_-$ of an
element $w \in \C(k)$ must be, and then define the positive tree
$T_+$ so that  the $w$ is given by the reduced tree pair
$ (T_-,T_+)$.  

\begin{figure}\includegraphics[width=3in]{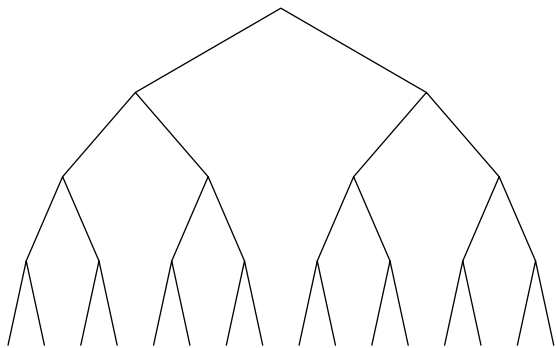}\\
\caption{The balanced tree $T_4$ \label{bal4}}
\end{figure}

Let $T_k$ be the balanced rooted binary tree with $2^k$ leaves;
that is, the tree with every node on the first $k$ levels having two
children, as in Figure \ref{bal4}.

For $w = (T_-,T_+)$ in the family $\C(k)$, we define $T_-$ to be the 
tree $T_{4k}$.  
Note that this is a very bushy tree, and has
at least $2k$ carets on the left side.  Each of these left carets has a
right subtree which is a complete tree with at least $k+2$ levels.  
Similarly,  $T_-$ has at least $2k$ right carets, each of which has a left 
subtree which is a complete tree with at least $k+2$ levels. 
There are a total of $2^{4k}$ leaves.

We construct the positive tree $T_+$ to have almost all carets
of type $\Ll$ and $\rni$, paired in a particular way with with
carets of $T_-$.  Let $r=2^{k-1}+2^{k-2}-1$ be the caret number of the first caret
on the right side of $T_-$.   Now let the tree $T_+$ correspond to the
word $x_0^{r-2} x_1 x_s$, where $s$ is $2^{4k}-3$.
Then $T_+$ will have  $2^{4k}$ leaves, the same number as in $T_-$.

\begin{figure}\includegraphics{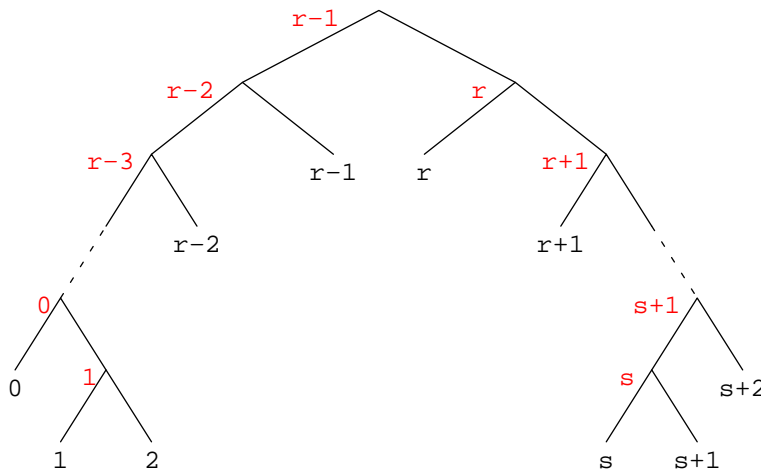}\\
\caption{Positive tree for a word $w \in \C(k)$ \label{posword}}
\end{figure}

We now check that with these definitions, $(T_-,T_+)$ forms a reduced 
tree pair diagram.  As pictured in Figure \ref{posword}, there are 
only  two carets in $T_+$ with two leaves: one with
leaves numbered $1$ and $2$, and the other with leaves numbered
$s=2^{4k}-3$ and $s+1=2^{4k}-2$.  In $T_-$, it is easy to see that 
because it is a complete tree, caret number $0$ has leaves numbered 
$0$ and $1$.  Also in $T_-$, the highest numbered caret has leaves 
numbered $s+1$ and $s+2$. Thus, no reduction of carets occurs, and $(T_-,T_+)$ is a reduced tree pair diagram.

In \S \ref{sec:action} above, the action of the generators of $F$
on a generic element is discussed.  We now describe the action of
a generator on an element $w$  of $\C(k)$, and more generally,
the action
of a sequence of generators on $w$.
 Let $\eta$ be a  word in the generators of $F$ which has
length strictly less than $k$, and let $w = (T_-,T_+)  \in \C(k)$.  We want to make
sure that $w \eta$ still satisfies the conditions of lemma
\ref{lemma:conditions}.  Because $\eta$ is not longer than $k$,
it can only affect a limited number of carets near the root of $T_-$.
For example, if $\eta$ is a power of $x_1$, then each application
of $x_1$ will rotate at  the right child $C_R$  of the root.
The left subtree of $C_R$ is, by construction, a 
complete tree with at least $k+2$ levels. Thus, after
performing $k$ clockwise rotations at the right child of the root, the resulting tree
still satisfies the conditions of lemma \ref{lemma:conditions}. 

More generally, no matter what the sequence of generators in $\eta$ is, the composition
of rotations that $\eta$ performs on $T_-$ affects carets only within distance
$k$ of the root.  Because of the fullness of the subtrees near 
the root of $T_-$,  the resulting tree will still have carets in the appropriate locations
to satisfy the conditions of  lemma \ref{lemma:conditions}. Because the exposed carets in $T_+$ are so far away from the root, we
know that no reductions can happen during the course of applying $\eta$ to $w$.
Thus,  lemma \ref{lemma:conditions}
guarantees that only one caret is affected by each application of the generator,

In the following chart we summarize the possible change in word
length when a generator of $\F$ acts on an element $w \eta$ with
$|\eta| < k$ and $w \in \C(k)$.  The positive tree $T_+$ has been chosen 
carefully
so that a caret in $w \eta$ affected by a generator is paired with
one of only two possible types of carets in $T_+$, an $\Ll$ or an
$\rni$.

\begin{center}

\begin{tabular}{|c|c|c|c|c|c|c|}

\hline Generator & Original & New  & Change in word  & Change in word\\
 & caret & caret & length when paired & length when paired \\
 & type  &  type& with $\Ll $ & with $\rni$\\ \hline
$\x$ & $\Ll$ & $R_I$ & -1 & \ 1 \\
\hline
$\ \inv$ & $R_I$ & $\Ll$ & \ 1 & -1  \\
\hline $x_1$ & $I_R$ & $R_I$ & -1 & -1 \\
\hline $\ x_1^{-1}$ & $R_I$ & $I_R$ &  \ 1 &  \ 1 \\
 \hline
\end{tabular}

\end{center}

We see immediately from this chart that $\x$ and $\inv$ will
reduce the word length of $w \in \C(k)$ because of the caret pairings
in $w$.  It is also true from the
chart that $x_1$ will reduce the length of the original word $w$.
The two elements we will consider to contradict almost convexity
will be $w \x$ and $w \inv$ for $w \in \C(k)$.  If the length  $|w|=n+1$,
then the length of $w \inv$ and $w \x$ will each be $n$.  Furthermore,
those two elements are distance 2 apart since there is an obvious path
from $w \x$ to $w$ to $w \inv$ of length 2.  That path, however, does
not lie in the ball of radius $n$.  In the proof of
theorem \ref{thm:notAC}, we will show that there is no short path from
$w \x $ to $w \inv$ which lies in  the ball of radius $n$.

\section{Almost convexity and $F$}

\label{sec:proofs}

We now prove that $F$ does not satisfy Cannon's $AC(2)$ condition,
and obtain as a corollary that $F$ does not satisfy $AC(n)$ for
any integral $n \geq 2$.

The idea of the proof of theorem \ref{thm:notAC} is the following.
Assuming $F$ satisfies the $AC(2)$ condition, we  would obtain a constant
$k$ so that any two points in $B(n)$ at distance $2$ from each
other would be connected by a path of length at most $k$ which
remains in $B(n)$.  Using this constant $k$, consider a point $w =
(T_-,T_+) \in \C(k+2)$.  The points $w \x$ and $w \inv$ are both in
$B(n)$ for $n= |w| - 1$ and are distance two apart.  Thus, there
would be a path $\g$ of length at most $k$ connecting them.  We assume this
path is oriented to go from $w \x$ to $w \inv$ and we follow the position
of the root caret $R$ of $T_-$ as it moves under the letters in
the path $\g$. We know that in $w \x$ the caret $R$ has moved to
the right side of the new negative tree.  The main lemma to the
proof of this theorem says that if at any time along the path $\g$
the caret $R$ becomes a left or an interior caret, then the path
$\g$ leaves $B(n)$ at that point.

Let $\g' = \x \g \x$ denote the loop based at $w$. The
contradiction to almost convexity arises from the following:
 Since the word  $w x_0$ has
$R$ as the right child of the root,  and the word $w x_0^{-1}$ has
$R$ as the left child of the root,  the final $\x$ in the path
$\g'$ would return $R$ to the root position from the left.  
Thus, at some point
along $\g$, the caret $R$ would have changed from
a right caret to a left or interior caret and  at that
point, the path $\g$ would have left the ball $B(n)$.

We begin with the proof of the necessary lemma.

\begin{lemma}

\label{lemma:outoftheball} Let $w=(T_-,T_+)  \in \C(k)$ with $|w| = n+1$, and $\g' =
\x^m \g'' \x$ be a loop based at $w$ of length at most $k$, with $m$ maximal.
 Let $R$ be the root caret of $T_-$, and $\eta$ the
shortest prefix of $\g''$ so that in $w \x^m \eta$ the caret $R$ is
not a right caret. Then the element $w \x^m \eta$ is not in
$B(n)$.
\end{lemma}

\begin{proof}

First, note that the negative tree of the element $w \x^m$ has
exactly $m$ right carets which are paired with $\Ll$ carets, and
we can number them as  we move away from the root as
$c_1,c_2, \cdots, c_m = R$, with
$c_1 < c_2 < \cdots < c_m$. Since the numbering
of the carets does not change when generators are applied, at the
first point where $R$ is not a right caret, then neither are any
of the carets $c_i$.

In the statement of the lemma, we are not distinguishing between
$R$ becoming a left caret and $R$ becoming an interior caret.  This will
not matter either for this proof or for the proof of theorem
\ref{thm:notAC} below.

The idea of the proof is to follow the path of each caret $c_i$ as
it is affected by different letters in the word $\eta$, and note
the net change in word length. Note that when we apply a generator
of $F$ to a word of the form $w \chi$, where $\chi$ is a word in
the generators of $F$ of length at most $k$, only a single caret
in the negative tree of $w \chi$ is affected.  In general, there
are times when this action can also affect a caret in the positive
tree, but we have chosen the form of elements of $\C(k)$ carefully
so that this is not the case, when applying strings of generators
of length less than $k$.

Each caret $c_i$ is originally paired with an $\Ll$ caret in the positive
tree by construction,
and since the positive tree will be unchanged, the positive part of
these pairing types will not change.
Consider all the letters in $\eta$ which change the caret type of
$c_i$. The last of these letters is either an $\inv$ changing
$c_i$ from a right caret to a left caret, or an $x_1^{-1}$
changing $c_i$ from a right caret to an interior caret.  According
to the chart in \S \ref{sec:family}, this is a net change in word
length of $+1$.

There are other letters in $\eta$ which can affect the caret
$c_i$. However, they must come in pairs, each pair leaving $c_i$
as a right caret so that the final letter in $\eta$ which affects
it can change it to a left or interior caret. These pairs can be
in one of two forms:

\begin{enumerate}

\item an $\inv$ which makes $c_i$ a left caret followed later in $\eta$
by an $\x$ making it again a right caret, or

\item an $x_1^{-1}$ making $c_i$ an interior caret and an $x_1$ later in $\eta$ 
 making it again a right caret.

\end{enumerate}

In either case, $c_i$ is always paired with an $\Ll$ caret, and we
see from the chart in \S \ref{sec:family} that the net change to
the total word length corresponding to either of these pairs is
always $0$. Thus, as we consider the letters of $\eta$ which
change the caret type of all $m$ of  the $c_i$'s, we see that they
contribute a total of $+m$ to the overall change in word length.

There may be letters in $\eta$ which affect the types of carets
other than the $c_i$.  Suppose caret $d \neq c_i$ is a caret affected by
a letter in $\eta$.  We claim that we  must have $d <
c_i$ for some $i$, and thus $d$ is also paired with an $\Ll$
caret.  If $d > c_i$ for all $i$, then $d$ would be  a caret which
appears after $R$.  In order for  $\eta$ to affect a caret after $R$,
 the caret $R$ would have had to have already moved from a right caret to a left or interior
caret, contradicting our assumption about $\eta$.  Thus,  we have established
the claim that  $d < c_i$ for some $i$.

Given the initial form of $w \in \C(k)$, we see that $d$ may begin
as an interior caret, and be initially moved to a right caret by
an element $x_1$.  From the chart in \S \ref{sec:family} we see
that this changes word length by $-1$.  Since $d<c_i$ for at least
one value of $i$, and all the $c_i$ must be changed from right
carets to non-right carets by the end of the path $\eta$, we must
also have $d$ changed from a right caret to a non-right caret.
Thus the last letter in $\eta$ affecting $d$ is either an $\inv$
which changes $d$ to a left caret or an $x_1^{-1}$ which changes
$d$ back to an interior caret. From the chart in \S
\ref{sec:family} we see that in either case, the change to the
word length is $+1$ making the total contribution of these two
letters in $\eta$ zero.

There may be other letters in $\eta$ which affect the caret $d$.
They must form the same pairs as listed above of ``intermediate"
letters which can affect the $c_i$, and thus contribute a total
word length change of zero.

The only other possibility for $d$ is that it begins as a left
caret, paired with an $\Ll$ caret for the same reasons as above.
Then the initial letter in $\eta$ affecting $d$ must be an $\x$,
making it a right caret.  The final letter in $\eta$
affecting $d$ again is either an $\inv$ or an $x_1^{-1}$.  Again,
we see from the chart in \S \ref{sec:family} that the net change
in word length coming from these two elements is $0$.  There can
also again be intermediate pairs of elements affecting $d$ of the same
forms as given above, which also contribute $0$ to the net change
in word length.

Since every letter of $\eta$ affects a single caret of $w$, each
letter of $\eta$ is one of the types listed above.  So the total
change in word length from $w\x$ to $w\x^m \eta$ is $m$.  Given
the initial form of $w$, it is easy to see from the chart that 
$|w \x^m| = |w| - m$.  Thus $|w \x^m \eta| = |w| - m + m = |w| = n+1$
and $ w \x^m \eta$ is not in $B(n)$.
\end{proof}

We are now ready to prove theorem \ref{thm:notAC} using the lemma
\ref{lemma:outoftheball}

\noindent {\it Proof of theorem \ref{thm:notAC}.} Assume that $F$
satisfies the $AC(2)$ condition.  Then there would be a constant $k$ so
that for every two points $x,y \in B(n)$ with $d(x,y) = 2$, there would
a be a path between them of length at most
$k$ lying completely inside $B(n)$.

Consider a point $w=(T_-,T_+) \in \C(k+2)$ with $|w| = n+1$.  By
construction, $|w \x| = |w \inv| = n$ and $d(w \x, w\inv) = 2$.
The assumption of almost convexity guarantees a path $\g$ from $w
\x$ to $w \inv$ whose length is bounded by $k$.  Let $\g' = \x \g
\x$ be the loop based at $w$ containing the path $\g$.

Let $R$ be the root caret in $T_-$.  The word  $w x_0$ has
$R$ as the right child of the root, so the initial $\x$ in the path
$\g'$ moves $R$ to a right caret.   The word $w x_0^{-1}$ has
$R$ as the left child of the root, so the final $\x$ in the path
$\g'$ must return $R$ to the root position from the left.  
Thus, at some point along the path $\g$, the caret
$R$ must change from being a right caret to a left caret.  So
there is a minimal prefix $\eta$ of $\g$ so that in $w \x \eta$,
the caret $R$ is not a right caret.  It then follows from lemma
\ref{lemma:outoftheball} that $w \x \eta$ is not in $B(n)$,
contradicting the assumption that $F$ is $AC(2)$. \qed

We immediately obtain the proof of corollary \ref{cor:notACn}.

\bibliography{thompc}
\bibliographystyle{plain}

\begin{small}
\noindent Sean Cleary \\
Department of Mathematics \\
City College of New York \\
City University of New York \\
New York, NY 10031 \\
E-mail: cleary@sci.ccny.cuny.edu

\medskip

\noindent
Jennifer Taback\\
Department of Mathematics and Statistics\\
University of Albany\\
Albany, NY 12222\\
E-mail: jtaback@math.albany.edu
\end{small}

\end{document}